\documentclass[final,leqno,onefignum,onetabnum]{siamltex1213}

\usepackage{amsmath,amssymb,apropos,pict2e,siunitx}

\title{Benchmark Computations of stresses in a spherical dome with shell finite elements\thanks{This work was
supported by the Academy of Finland through the project \emph{How to handle the
prima donna of structures: analysis and development of advanced discretization techniques
for the simulation of thin shells} (decision number 260302).}}

\author{Antti H.~Niemi\thanks{Aalto University, Espoo, Finland. 
(\email{antti.h.niemi@aalto.fi}). Questions, comments, or corrections
to this document may be directed to that email address.}}

\begin{document}
\maketitle
\slugger{sic}{xxxx}{xx}{x}{x--x}

\begin{abstract}
We present a computational framework for analysing thin shell structures using the finite element method. The framework is based on a mesh-dependent shell model which we derive from the general laws of three-dimensional elasticity. We apply the framework for the so called Girkmann benchmark problem involving a spherical shell stiffened with a foot ring. In particular, we compare the accuracy of different reduced strain four-node elements in this context. We conclude that the performance of the bilinear shell finite elements depends on the mesh quality but reasonable accuracy of the quantities of interest of the Girkmann problem can be attained in contrast to earlier results obtained with general shell elements for the problem.
\end{abstract}



\pagestyle{myheadings}
\thispagestyle{plain}
\markboth{ANTTI H. NIEMI}{BENCHMARK COMPUTATIONS OF STRESSES IN A SPHERICAL DOME}

\section{Introduction}
The role of thin shell theory in structural analysis has changed dramatically over centuries from center stage to supporting cast, partly because of the advent of the finite element method. This paper will bring shell theory back to the forefront by studying a family of four-node finite elements based on it. The study is carried out in the context of a challenging benchmark problem involving a stiffened doubly curved thin shell.

Conventional shell finite element formulations involve various explicit and implicit modelling assumptions that extend beyond the limits of mathematical convergence theory currently available. The theoretical problems arise from the fact that the general shell elements used in industrial FEA have been developed through ``finite element modelling'' where the connection of the discretization to the actual differential operators of the mathematical shell model is obscure, cf.~\cite{Szabo2012}.

A prerequisite for traditional FE error analysis is a well-posed variational problem formulated in some suitable Hilbert space. For shells such mathematical models can be formulated using differential geometry of surfaces and theoretically stable formulations have been analyzed e.g.~in \cite{Pitkaranta1992,Arnold1997,Chapelle1998a}. However, these works deal only with the bending-dominated deformations and do not address the membrane-dominated and intermediate cases which are very important for instance in civil and structural engineering.

On the other, it is possible to interpret the modelling assumptions of conventional shell elements in context of the mathematical models. For instance, the degenerated solid approach associated to certain four-node elements has been translated to explicit strain reduction procedure within a specific shell model in \cite{Malinen2001,Malinen2002} and that procedure has been numerically analysed in \cite{Havu2001,Havu2003,Niemi2008a,Niemi2008b}. This line of research is not limited to the theoretical analysis of existing formulations only. It can also be used to enhance the formulations and develop new ones as shown in \cite{Niemi2010}, where a
new four-node shell finite element of arbitrary quadrilateral shape was developed based on shell theory.

Thanks to modern computation technology, such as the $hp$-adaptive finite element method \cite{demkowicz2006,tews2009}, shell analysis can also be based directly to three-dimensional elasticity theory. Such an approach rules out the modelling errors arising from the simplifications of dimensionally reduced structural models but requires more degrees of freedom for the discrete model. Also, if simplified representations of the stress state such as the stress resultants are needed, they must be post-processed from the three-dimensional stress field and this can be non-trivial.

A model problem called the \emph{Girkmann problem}, which was revived some time ago, highlights the above complications rather dramatically, see \cite{Szabo2009,Niemi2012,Pitkaranta2012,Devloo2013}. The problem involves a concrete structure consisting of a spherical dome stiffened by a foot ring under a dead gravity load. The task is to determine the values of the shear force and the bending moment at the junction between the dome and the ring as well as the maximum bending moment in the dome. 

The problem was initially presented and solved analytically in the text book \cite{Girkmann1963}. More recently, in the bulletin of the International Association of Computational Mechanics (IACM) \cite{Pitkaranta2008}, the problem was posed as a computational challenge to the finite element community. The purpose of the challenge was to find out how the process of verification, that is the process of building confidence that an approximative result is within a given tolerance of the exact solution to the mathematical model, is carried out by the community given the Girkmann problem. The results, that are summarized in \cite{Pitkaranta2009,Szabo2009}, without attribution and details on how verification was actually performed, are scandalous. Out of the 15 results submitted, 11 have a very large dispersion and are not within any acceptable tolerance of the reference values computed in \cite{Szabo2009,Niemi2012,Pitkaranta2012} using different models and formulations.

So far detailed verification studies have been published for the axisymmetric models based on elasticity theory as well as axisymmetric dimensionally reduced models. In \cite{Szabo2009}, the $p$-version of FEM was used in conjunction with the extraction procedure of \cite{Babuska1984} to compute accurate values for the quantities of interest. Similar approach with the $hp$-version of FEM was taken in \cite{Niemi2012}, where also the axisymmetric $h$-version with selective reduced integration was successfully employed to discretize the dimensionally reduced model.

In the present work, we introduce a finite element framework for thin shell analysis and 
use the Girkmann problem to demonstrate its possibilities. More precisely, we benchmark different variants of the MITC-type shell element proposed earlier by the author in \cite{Niemi2010} by modelling a quarter of the dome and by using symmetry boundary conditions. We show that the performance of the formulations varies depending on the performed strain reductions and the mesh regularity. Nevertheless, reasonable accuracy for the main quantities of interest is obtained contrary to the earlier published results obtained with general shell elements.

The paper is organized as follows. In the next Section, we develop the shell theory used to construct the different FE methods. The finite element methods are described in Section \ref{sec:FEM} together with a discussion on their well-posedness and implementation aspects. Section \ref{sec:numerics} is devoted to description of the analysis procedure for the Girkmann problem and numerical results. The paper ends with conclusions in Section \ref{sec:conclusions}.

\section{Variational Formulation of the Shell Problem}
We will employ the Einstein summation convention so that Greek indices range over the values 1, 2 while Latin indices have the values 1, 2, 3. The former refer to \emph{surface coordinates} while the latter refer to general three-dimensional curvilinear coordinates. Our surface coordinate systems can be assumed orthonormal so that we can formulate our strain-displacement relations and constitutive laws directly in terms of physical components and traditional partial derivatives which will be denoted by a comma. Moreover, Euclidean vectors are displayed using overhead arrows, whereas boldface notation is reserved for the column vectors and matrices  storing the components of different surface tensors in the assumed orthonormal coordinate system.

A shell domain $\Omega \subset \Reals^3$ of constant thickness $t$ is defined as
\begin{equation} \label{eqn:shell domain}
\Omega = \mPhi(K \times (-t/2,t/2)),
\end{equation}
where the shell mapping $\mPhi$ is of the form
\begin{equation} \label{eqn:shell mapping}
\mPhi(x,y,\zeta) = \vr(x,y) + \zeta \vn(x,y).
\end{equation}
Here the parametric surface $\vr(x,y)$ represents the middle surface of the shell and $\vn(x,y)$ is the unit normal vector to the middle surface. Thus $\xi_1 = x$, $\xi_2 = y$, and $\xi_3 = \zeta$ constitute a curvilinear coordinate system in 3-space, called \emph{shell coordinates}.

In the following, we imagine that the middle surface is defined as
\begin{equation} \label{eqn:middle surface parametrization}
\vr(x,y) = x \vec i_1 + y \vec i_2 + f(x,y)\vec i_3, \quad (x,y) \in K,
\end{equation}
where $\vec i_1,\vec i_2,\vec i_3$ are \emph{fixed} Cartesian unit vectors. Moreover, we assume that the middle surface differs only little from the coordinate plane $K$, i.e.~the shell is \emph{shallow}. More precisely, we assume that the function $f:K \rightarrow \Reals$ is smooth in the curvature length scale $R$ defined by
\[
\frac{1}{R}=\max_{\alpha,\beta} \norm{f_{,\alpha\beta}}_{\infty,K}.
\] 
If $h_K=\mathrm{diam}(K)$ and $R$ is taken as the length unit, the shallowness assumption can be formulated as
\begin{equation} \label{eqn:shallowness assumption}
f_{,\alpha} = \cO(\hat{h}_K), 
\end{equation}
where $\hat{h}_K=h_K/R \leq 1$. 

Under the shallowness assumption \eqref{eqn:shallowness assumption}, the tangent basis vectors 
\begin{equation} \label{eqn:tangent vectors}
\ve_\alpha(x,y) = \vec i_\alpha + f_{,\alpha}(x,y) \vec i_3
\end{equation}
associated to the middle surface parametrization \eqref{eqn:middle surface parametrization}, are orthonormal within the accuracy of $\cO(\hat{h}_K^2)$.

\subsection{Shell Kinematics}
According to the standard kinematic hypothesis we assume that the displacement vector can be written in the form
\begin{equation} \label{eqn:displacement hypothesis}
\vU(x,y,\zeta) = (u_\lambda(x,y)+\zeta\theta_\lambda(x,y))\ve_\lambda(x,y) + w(x,y)\vn(x,y),
\end{equation}
where $\bu = (u_1,u_2)$ are the tangential displacements of the middle surface, $w$ is the transverse deflection, and the quantities $\bstheta = (\theta_1,\theta_2)$ are the angles of rotation of the normal.

Referring to the curvilinear coordinates $\xi_1 = x$, $\xi_2 = y$, and $\xi_3 = \zeta$, the linearized Green-Lagrange strain tensor is defined by
\begin{equation} \label{eqn:3D strain tensor}
e_{ij} = \frac{1}{2}(\mPhi_{,i} \cdot \vU_{,j} + \mPhi_{,j} \cdot \vU_{,i}), \quad i,j=1,2,3.
\end{equation}
We find directly from \eqref{eqn:displacement hypothesis} that
\begin{equation} \label{eqn:displacement tangential derivatives}
\vU_{,\alpha} = (u_{\lambda,\alpha} + \xi_3 \theta_{\lambda,\alpha})\ve_\lambda + (u_\lambda + \xi_3\theta_\lambda)\ve_{\lambda,\alpha} + w_{,\alpha} \vn + w \vn_{,\alpha}
\end{equation}
for $\alpha=1,2$ and
\begin{equation} \label{eqn:displacement normal derivative}
\vU_{,3} = \theta_\lambda \ve_\lambda.
\end{equation}
Similarly, combination of \eqref{eqn:shell mapping}, \eqref{eqn:middle surface parametrization} and \eqref{eqn:tangent vectors} yields
\begin{equation} \label{eqn:shell mapping tangential derivatives}
\mPhi_{,\alpha} = \ve_{\alpha} + \zeta \vn_{,\alpha}, \quad \alpha=1,2,
\end{equation}
and 
\begin{equation} \label{eqn:shell mapping normal derivative}
\mPhi_{,3} = \vn.
\end{equation}

The components of the strain tensor are represented as power series of the variable $\zeta$. If we take into account two terms, the in-plane strains may be written as
\begin{equation} \label{eqn:in-plane strains}
e_{\alpha\beta} \approx \varepsilon_{\alpha\beta} + \zeta \kappa_{\alpha\beta}.
\end{equation}
Using relations \eqref{eqn:shell mapping tangential derivatives} and \eqref{eqn:displacement tangential derivatives}, the membrane strain tensor $\varepsilon_{\alpha\beta}$, which arises from stretching of the deformed middle surface, can be written as
\begin{equation} \label{eqn:membrane strain tensor}
\varepsilon_{\alpha\beta} \approx \frac{1}{2}(u_{\alpha,\beta}+u_{\beta,\alpha}) - b_{\alpha\beta} w,
\end{equation}
where 
\begin{equation} \label{eqn:curvature tensor}
b_{\alpha\beta} = - \ve_\alpha \cdot \vn_{,\beta} = \frac{f_{,\alpha\beta}}{\sqrt{1+f_{,\lambda}f_{,\lambda}}}, \quad \alpha,\beta=1,2,
\end{equation} 
are the coefficients of the second fundamental form of the middle surface. In \eqref{eqn:membrane strain tensor}, terms multiplied by the `rotation coefficients' $\ve_\alpha \cdot \ve_{\lambda,\beta}$ have been neglected as quantities of relative order $\cO(\hat{h}_K)$ based on the shallowness assumption \eqref{eqn:shallowness assumption}. We content ourselves here to first order of accuracy since, in general, the coefficients  $b_{\alpha\beta}$ cannot be approximated more precisely with linear or bilinear interpolation functions, see Section~\ref{subsec:implementation}.

Introducing the coefficients of the third fundamental form of the middle surface
\[
c_{\alpha\beta} = \vn_{,\alpha} \cdot \vn_{,\beta}, \quad \alpha,\beta=1,2,
\]
the elastic curvature tensor $\kappa_{\alpha\beta}$, which arises from bending of the deformed middle surface, comes out as
\begin{equation} \label{eqn:bending strain tensor}
\kappa_{\alpha\beta} \approx \frac{1}{2}(\theta_{\alpha,\beta} + \theta_{\beta,\alpha}) + c_{\alpha\beta} w - \frac{1}{2} (b_{\alpha\lambda}u_{\lambda,\beta} + b_{\beta\lambda}u_{\lambda,\alpha}), \quad \alpha,\beta=1,2
\end{equation}
based again directly on \eqref{eqn:shell mapping tangential derivatives} and \eqref{eqn:displacement tangential derivatives}. Here, terms multiplied by $\ve_\alpha \cdot \ve_{\lambda,\beta}$ or $\vn_{,\alpha} \cdot \ve_{\lambda,\beta}$ have been neglected as quantities of relative order $\cO(\hat{h}_K)$.

It is possible to simplify the bending strain expressions by sacrificing their tensorial invariance. It is straightforward to verify that
\[
c_{\alpha\beta} \approx b_{\alpha\lambda} b_{\lambda\beta}
\]
within the adopted accuracy, so that we may write \eqref{eqn:bending strain tensor} component-wise as
\[
\begin{aligned}
\kappa_{11} &\approx \underline{\theta_{1,1} + b_{12} (b_{12}w - u_{2,1}}) - b_{11} \varepsilon_{11}, \\
\kappa_{22} &\approx \underline{\theta_{2,2} + b_{12}(b_{12} w - u_{1,2})} -b_{22} \varepsilon_{22}, \\
\kappa_{12} &\approx \underline{\frac{1}{2}\left( \theta_{1,2} + \theta_{2,1} 
+ b_{11} (b_{12} w - u_{1,2}) + b_{22} (b_{12} w -u_{2,1}) 
\right)}  -\frac{ b_{12}}{2} (\varepsilon_{11} + \varepsilon_{22}). \\
\end{aligned} 
\]
In these expressions, the contribution of the terms $b_{11} \varepsilon_{11}$, $b_{22} \varepsilon_{22}$ and $b_{12}(\varepsilon_{11}+\varepsilon_{22})$ to the maximum in-plane strains at the outer and inner surfaces of the shell is of relative order $\cO(t/R)$ only. Therefore, the number of terms in the kinematic relations can be slightly reduced by retaining only the underlined terms in the calculations.

Finally, the transverse shear strains are defined as 
\begin{equation} \label{eqn:out-of-plane strains}
\gamma_\alpha = 2e_{\alpha3}, \quad \alpha = 1,2.
\end{equation}
These can be written in terms of the displacements by first noting that since $\ve_\alpha$ and $\vn$ are orthogonal, we have $\ve_{\alpha,\beta} \cdot \vn = -\ve_{\alpha} \cdot \vn_{,\beta}$, and consequently $b_{\alpha\beta} = \vn \cdot \ve_{\alpha,\beta}$. Now, combination of \eqref{eqn:shell mapping normal derivative} with \eqref{eqn:displacement tangential derivatives} and \eqref{eqn:shell mapping tangential derivatives} with \eqref{eqn:displacement normal derivative} according to \eqref{eqn:3D strain tensor} yields 
\begin{equation} \label{eqn:transverse shear strains}
\gamma_\alpha = \theta_\alpha + b_{\lambda\alpha} u_\lambda + w_{,\alpha}, \quad \alpha=1,2,
\end{equation}
and completes the description of shell kinematics.

\subsection{Constitutive Equations}
Assuming linearly elastic isotropic material with Poisson ratio $\nu$ and Young modulus $E$, the constitutive law relating stresses and strains can be written with respect to the approximately orthogonal shell coordinate system $(x,y,\zeta)$ attached to the middle surface as
\begin{equation} \label{eqn:constitutive law}
\begin{aligned}
\sigma_{\alpha\beta} &= \frac{E}{1-\nu^2} \left[ (1-\nu) e_{\alpha\beta} + \nu e_{\lambda\lambda}\delta_{\alpha\beta} \right], \\
\sigma_{\alpha3} &= \frac{E}{1+\nu}e_{\alpha3},\quad \alpha,\beta = 1,2,
\end{aligned}
\end{equation}
where $\delta_{\alpha\beta}$ is the Kronecker delta. The elastic coefficients in the above formula have been modified to yield so called plane stress state tangent to the middle surface. This modification is necessary to avoid Poisson locking in context of the kinematic assumption \eqref{eqn:displacement hypothesis}.

We follow the standard convention of structural mechanics and introduce the internal forces and moments which are the stress resultants and stress couples per unit length of the middle surface. These can now be defined as
\begin{equation} \label{eqn:stress resultants}
\begin{aligned}
n_{\alpha\beta} = \int_{-t/2}^{t/2} \sigma_{\alpha\beta}\,d\zeta, \quad
m_{\alpha\beta} = \int_{-t/2}^{t/2} \sigma_{\alpha\beta}\zeta\,d\zeta, \quad
q_\alpha = \int_{-t/2}^{t/2} \sigma_{\alpha 3}\,d\zeta,
\end{aligned}
\end{equation}
and correspond to the membrane forces, bending moments and transverse shear forces in static equilibrium considerations.

\subsection{Potential Energy Functional}
Integrals over the shell domain $\Omega$ can be evaluated in terms of the assumed shell coordinates as
\[
\int_{\Omega} \{\cdot\}\,d\Omega \approx \int_K \int_{-t/2}^{t/2} \{\cdot\}\,d\zeta dx dy.
\]
Combination of \eqref{eqn:stress resultants}, \eqref{eqn:constitutive law}, \eqref{eqn:out-of-plane strains} and \eqref{eqn:in-plane strains} allows us to write the elastic strain energy functional as
\begin{equation} \label{eqn:strain energy}
\begin{split}
U_K(\bu,w,\bstheta) = \frac{1}{2} \int_{\Omega} \sigma_{\alpha\beta} e_{\alpha\beta} \,d\Omega 
\approx \frac{1}{2} &\int_K (n_{\alpha\beta} \varepsilon_{\alpha\beta}
+ q_\alpha \gamma_\alpha
+ m_{\alpha\beta}\kappa_{\alpha\beta})\,dx dy,
\end{split}
\end{equation}
where 
\begin{equation} \label{eqn:explicit stress resultants}
\begin{aligned}
n_{\alpha\beta} &= \frac{Et}{1-\nu^2} \left[ (1-\nu) \varepsilon_{\alpha\beta} + \nu \varepsilon_{\lambda\lambda}\delta_{\alpha\beta} \right], \\
q_{\alpha} &= \frac{Et}{2(1+\nu)} \gamma_\alpha, \\
m_{\alpha\beta} &= \frac{Et^3}{12(1-\nu^2)} \left[ (1-\nu) \kappa_{\alpha\beta} + \nu \kappa_{\lambda\lambda}\delta_{\alpha\beta} \right]
\end{aligned}
\end{equation}
and the strains are given in terms of the displacements in \eqref{eqn:membrane strain tensor}, \eqref{eqn:transverse shear strains} and \eqref{eqn:bending strain tensor}.

Similarly, the potential energy corresponding to external distributed surface forces $(f_1,f_2,p)$ and moments $(\tau_1,\tau_2)$ is 
\begin{equation}
V_K(\bu,w,\bstheta) = -\int_K (f_\lambda u_\lambda + pw + \tau_\lambda \theta_\lambda)\,dx dy
\end{equation}
and the total energy is given by the sum
\begin{equation}
E_K(\bu,w,\bstheta) = U_K(\bu,w,\bstheta) + V_K(\bu,w,\bstheta).
\end{equation}

\section{Finite Element Methods} \label{sec:FEM}
We assume that the whole shell domain $\tilde{\Omega}$ is formed as a union of domains of the form \eqref{eqn:shell domain} as
\[
\tilde{\Omega} = \bigcup_{K \in \cC_h} \Omega_K,
\]
where $\cC_h$ stands for a mesh of convex quadrilaterals $K$ corresponding to parametrizations of patches of the shell middle surface according to \eqref{eqn:middle surface parametrization}. We also assume that the whole middle surface
\[
S = \bigcup_{K \in \cC_h} \vr_K(K)
\] 
is a smooth surface and that it can be described alternatively by a single, global parametrization $\vrho(\tilde{\xi}_1,\tilde{\xi}_2)$. It follows that the transformations between the local and global coordinate systems $T_K = \vr_K^{\,-1} \circ \vrho$, $K \in \cC_h$, are diffeomorphisms.

Without losing generality, we may assume that the coordinates $\tilde{\xi}_1,\tilde{\xi}_2$ are isothermal. If $\{\vg_1(\tilde{\xi}_1,\tilde{\xi}_2),\vg_2(\tilde{\xi}_1,\tilde{\xi}_2) \}$ are the corresponding orthonormal tangent vectors of the middle surface, and $\tilde{u}_\alpha$, $\tilde{w}$, and $\tilde{\theta}_\alpha$ stand for the associated displacement and rotation components, then these components are related to the local components $u_\alpha^K$, $w^K$ and $\theta_\alpha^K$ as
\begin{equation} \label{eqn:displacement transformation}
u_\alpha^K \circ T_K = \tilde{u}_\lambda \vg_\lambda\cdot \vi_\alpha, \quad 
w^K \circ T_K = \tilde{w}, \quad \theta_\alpha^K \circ T_K = \tilde{\theta}_\lambda \vg_\lambda \cdot \vi_\alpha
\end{equation}
according to \eqref{eqn:tangent vectors} and \eqref{eqn:displacement hypothesis}. The total potential energy of the structure is expressed as the sum of element-wise contributions such that
\begin{equation}
E(\tilde{\bu},\tilde{w},\tilde{\bstheta}) = \sum_{K \in \cC_h} E_K(\bu^K,w^K,\bstheta^K),
\end{equation}
where $(\tilde{\bu},\tilde{w},\tilde{\bstheta})$ is the globally defined generalized displacement field. The solution of the problem is determined according to the principle of minimum potential energy from the condition
\[
E(\tilde{\bu},\tilde{w},\tilde{\bstheta}) = \min_{(\tilde{\bv},\tilde{z},\tilde{\bspsi}) \in \cU} E(\tilde{\bv},\tilde{z},\tilde{\bspsi}),
\]
where the energy space $\cU$ is defined as the set of those kinematically admissible generalized displacement fields for which the energy functional is finite. The existence of a unique minimizer (for $\max_K h_K$ sufficiently small) follows from the well-posedness of the corresponding Reissner-Naghdi shell model, see e.g.~\cite{Chapelle2011}. 

It is now straightforward to formulate a finite element method where each displacement component is approximated separately as in the space
\begin{equation}
\cU_h = \{ (\tilde{\bu}, \tilde{w}, \tilde{\bstheta}) \in \cU \; : \; 
(\bu^K, w^K, \bstheta^K) \in [Q_1(K)]^5 \; \forall K \in \cC_h \},
\end{equation}
where $Q_1(K)$ denotes the standard space of isoparametric bilinear functions on $K$ and $(\bu^K, w^K, \bstheta^K)$ is the local generalized displacement field defined in \eqref{eqn:displacement transformation}.

\subsection{Strain Reduction Techniques}
To avoid locking when approximating bending-dominated problems, membrane and transverse shear strains must be reduced. To introduce the different methods, we denote by $\bF_K = (x_K,y_K)$ the bilinear mapping of the reference square $\hat{K} = [-1,1] \times [-1,1]$ onto $K$ and by 
\[
\bJ_K = \bs{J}_K(\hat{x},\hat{y}) =  \begin{pmatrix}
\dfrac{\partial x_K}{\partial \hat{x}} & \dfrac{\partial x_K}{\partial \hat{y}} \\ 
\dfrac{\partial y_K}{\partial \hat{x}} & \dfrac{\partial y_K}{\partial \hat{y}} \\
\end{pmatrix}
\]
the Jacobian matrix of $\bF_K$. Here $(\hat{x},\hat{y})$ are the coordinates on $\hat{K}$.

We start by defining on the reference square $\hat{K}$ the function spaces 
\begin{equation} \label{eqn:reference shear strain space}
\bS(\hat{K}) = \{ \bs{\hat{s}} = \begin{pmatrix} a + b\hat{y} \\ c + d\hat{x} \end{pmatrix} \; : \; a,b,c,d \in \Reals \}
\end{equation} 
and 
\begin{equation} \label{eqn:reference membrane strain space}
\bM(\hat{K}) = \{ \boldsymbol{\hat{\tau}} = \begin{pmatrix} a + b\hat{y} & c \\ c & d + e\hat{x} \end{pmatrix} \; : \; a,b,c,d,e \in \Reals \}
\end{equation}
for the reduced transverse shear strains and membrane strains, respectively. The canonical degrees of freedom associated with $\bs{S}(\hat{K})$ are
\begin{equation} \label{eqn:shear strain dofs}
\bs{\hat{s}} \mapsto \int_{\hat{e}} \bs{\hat{s}}^T \bs{\hat{t}}\,d\hat{s} \; \text{for every edge $\hat{e}$ of $\hat{K}$},
\end{equation}
whereas the degrees of freedom associated with $\bs{M}(\hat{K})$ are defined as
\begin{equation} \label{eqn:membrane strain dofs}
\begin{aligned}
\bs{\hat{\tau}} &\mapsto \int_{\hat{e}} \bs{\hat{t}}^T \bs{\hat{\tau}}\bs{\hat{t}}\,d\hat{s}\; \text{for every edge $\hat{e}$ of $\hat{K}$}, \\ 
\bs{\hat{\tau}} &\mapsto \int_{\hat{K}} \hat{\tau}_{12}\,d\hat{x}d\hat{y}.
\end{aligned}
\end{equation}

The corresponding spaces associated to a general element $K \in \cC_h$ are then defined using covariant transformation formulas as
\begin{equation} \label{eqn:covariant shear transformation}
\bS(K) = \{ \bs{s} = (\bJ_K^{-T} \boldsymbol{\hat{s}}) \circ \bF_K^{-1} = \cS_K(\bs{\hat{s}}) \; : \; \boldsymbol{\hat{s}} \in \bS(\hat{K}) \}
\end{equation}
and
\begin{equation} \label{eqn:covariant membrane transformation}
\bM(K) = \{ \bs{\tau} = \bs{\bar{J}}_K^{-T} (\boldsymbol{\boldsymbol{\hat{\beta}}} \circ \bF_K^{-1})\bs{\bar{J}}_K^{-1} = \bar{\cM}_{K}(\bs{\hat{\tau}}) \; : \; \boldsymbol{\hat{\tau}} \in \bM(\hat{K}) \}
\end{equation}
In case of general quadrilateral elements with non-constant Jacobian matrices, the transformation in \eqref{eqn:covariant membrane transformation} must be fixed to a single orientation e.g.~at the midpoint of $\hat{K}$ so that $\bs{\bar{J}} = \bs{J}(0,0)$. Otherwise the function space \eqref{eqn:covariant membrane transformation} for the reduced membrane strains and stresses does not necessarily include constant fields which would degrade the accuracy of the formulation. 

Denoting by $\bs{\mathit{\Lambda}}_{\hat{K}}:\bs{H}^1(\hat{K}) \rightarrow \bs{S}(\hat{K})$ and $\bs{\mathit{\Pi}}_{\hat{K}}:\bs{H}^1(\hat{K}) \rightarrow \bs{M}(\hat{K})$ the interpolation operators associated to the degrees of freedom \eqref{eqn:shear strain dofs} and \eqref{eqn:membrane strain dofs}, the corresponding projectors for a general $K \in \cC_h$ are defined as
\[
\mPi_K = \bar{\cM}_K \circ \mPi_{\hat{K}} \circ \bar{\cM}_K^{-1} \; \text{and} \; \mLambda_K = \cS_K \circ \mLambda_{\hat{K}} \circ \cS_K^{-1}
\]
The transformation rule \eqref{eqn:covariant shear transformation} guarantees that the degrees of freedom \eqref{eqn:shear strain dofs} are preserved on $K$:
\[
\bs{s} \mapsto \int_e \bs{s}^T\bs{t}\,ds = 0 \; \text{for every edge $e$ of $K$},
\]
but the analogous statement does not hold in context of \eqref{eqn:covariant membrane transformation} for general meshes.

The finite element introduced in \eqref{eqn:reference shear strain space} and \eqref{eqn:shear strain dofs} is a well known edge element denoted by the symbol RTc$_1^e$ in the recently introduced Periodic Table of the Finite Elements \cite{Arnold2014} and has been described in the context of plate bending e.g.~in \cite{Hughes1981,Bathe1985}. The element \eqref{eqn:reference shear strain space}, \eqref{eqn:membrane strain dofs} is less customary at least in the mathematical literature -- probably because of its non-elegant extension to quadrilateral shapes. In any case the  historical roots of the element are very deep in the literature on finite element technology for plane elasticity. Our current formulation corresponds essentially to the stress field of the Pian-Sumihara element introduced in \cite{Pian1984} which in turn may be viewed as an extension of the nonconforming displacement methods introduced by Wilson et.~al.~\cite{Wilson1973} and Turner et.~al.~\cite{Turner1956}. We refer the reader to \cite{Pitkaranta2000a} for the complete mathematical theory of these formulations in context of plane elasticity. In the present context, the convergence theory is confined to special cases involving cylindrical or globally shallow shells on special meshes, see~\cite{Havu2001,Havu2003,Niemi2008a,Niemi2008b}.

We shall use the label MITC4C for the formulation for which only the transverse shear strains are projected into the space \eqref{eqn:covariant shear transformation} 
\begin{equation}
\bs{\gamma} \hookrightarrow \Lambda_K \bs{\gamma}
\end{equation}
and the label MITC4S for the formulation where also the membrane strains are projected:
\begin{equation}
\bs{\varepsilon} \hookrightarrow \mPi_K \bs{\varepsilon}, \quad \bs{\gamma} \hookrightarrow \Lambda_K \bs{\gamma}
\end{equation}
when evaluating the strain energy according to \eqref{eqn:strain energy} and \eqref{eqn:explicit stress resultants}. 

In addition, we consider stabilized variants of both methods, where the shear modulus $G = \frac{E}{2(1+\nu)}$ in \eqref{eqn:explicit stress resultants} is modified as
\begin{equation} \label{eqn:shear balancing}
G \hookrightarrow G_K = \frac{t^2}{t^2+\alpha_K h_K^2} \cdot G.
\end{equation}
Here $\alpha_K$ is a positive stabilization parameter independent of $t$ and $h_K$. This stabilization idea originates from the corresponding plate bending elements, see e.g.~\cite{Tessler1983,Lyly1993}.

Finally, the abbreviation DISP4 is used for the standard displacement method without any strain reduction or stabilization.

\subsection{Implementation} \label{subsec:implementation}
The implementation of the present formulation follows the standard steps used to construct quadrilateral plane-elastic elements with some twists. The shape functions are the usual bilinear ones on the reference square $\hat{K}$ and the computation of derivatives and numerical integration using the Gauss quadrature are performed in the canonical way. 

The first twist concerns determination of the physical element $K$ from the general surface mesh consisting of quadrilateral elements. As it may happen, that the four nodes of an element do not lie in the same plane, a straightening operation is needed in order to define the plane element $K$. This can be accomplished in many ways, but we follow the procedure of \cite{Macneal1978} described also in \cite{Niemi2010} that yields naturally also the directions of the local coordinate axes $\vec i_1$ and $\vec i_2$.

The second difference is related the strain-displacement relations \eqref{eqn:membrane strain tensor},\eqref{eqn:bending strain tensor} and \eqref{eqn:transverse shear strains}, which involve the coefficients of the second fundamental form in addition to the standard shape function derivatives. Within the limits of the local shallowness assumption, these coefficients can computed using the \emph{interpolated normal vector} $\vec n_h$ as
\[
b_{\alpha\beta} \approx - \vec i_\alpha \cdot \vec n_{h,\beta}, \quad \alpha,\beta=1,2.
\]
To carry out the interpolation, the nodal normals are needed as geometric input data in addition to the node positions. The former are also used to construct the two orthonormal tangent direction $\vec g_1,\vec g_2$ used to enforce the continuity of the tangential displacements and normal rotations according to \eqref{eqn:displacement transformation}.

The current implementation (assembly, solution and post-processing) is carried out using Mathematica while Gmsh is used for the mesh generation \cite{WolframResearch2015,Geuzaine2009}. An additional pre-processing step is the determination of the nodal normals that can be performed by using the analytic surface representation if available, or by averaging the normals of the elements sharing a common node.

\section{Numerical Results} \label{sec:numerics}
\begin{figure}
\begin{center}
\setlength{\unitlength}{0.75cm}
\begin{picture}(16,9)(0,3)
\thicklines
\put(0,-10){\arc[90,50]{20}}
\put(12.856,5.321){\line(1,0){0.514}}
\put(12.856,4.892){\line(1,0){0.514}}
\put(12.856,4.892){\line(0,1){0.429}}
\put(13.370,4.892){\line(0,1){0.429}}
\thinlines
\put(0,4.892){\vector(1,0){16}}
\put(0,4.892){\vector(0,1){6}}
\put(12.856,5.321){\line(0,1){3}}
\put(10.856,7.5){\vector(1,0){2}} \put(12.2,7.8){$\rho_0$}
\put(5.325,5.974){\vector(1,3){1}}
\put(6.35,8.2){$r_0$}
\put(13,11){$\rho_0=15.0$\ m}
\put(13.1,10){$\,\alpha=40^\circ$}
\put(13,9){$\,r_0=\rho_0 / \sin\alpha$}
\end{picture}
\setlength{\unitlength}{0.75cm}
\begin{picture}(15,6)(-1.5,0)
\thicklines
\put(2,4.164){\line(-4,3){2}}
\put(2.252,4.5){\line(-4,3){2}}
\put(2,4.164){\line(3,4){0.252}}
\put(8,1){\line(1,0){4.2}}
\put(8,1){\line(0,1){3.164}}
\put(8.252,4.5){\line(1,0){3.948}}
\put(12.2,1){\line(0,1){3.5}}
\put(8,4.164){\line(3,4){0.252}}
\thinlines
\put(9.5,0.35){$0.60\,\mathrm{m}$}
\put(12.5,2.5){$0.50\,\mathrm{m}$}
\put(7.85,0.5){$A$}
\put(12.15,0.5){$B$}
\put(2.126,4.332){\vector(4,-3){2}} \put(4.2,2.55){$N$}
\put(2.126,4.332){\vector(3,0){2}} \put(4.2,4.2){$R$}
\put(2.35,4.2){\arc[-140,40]{0.5}} \put(2.68,4.6){\vector(-1,1.5){0.03}}
\put(2.4,4.75){$M$}
\put(8.126,4.332){\vector(-4,3){2}} \put(5.7,5.75){$N$}
\put(8.126,4.332){\vector(-3,0){2}} \put(5.7,4.2){$R$}
\put(8,4.5){\arc[-120,-305]{0.5}} \put(8.4,4.85){\vector(1,-0.6){0.05}}
\put(8.5,4.75){$M$}
\put(0.4,2){\vector(3,4){1.5}}
\put(0.4,2){\line(0,1){2}}
\put(0.4,2){\arc[90,65]{1}}
\put(0.91,2.82){\vector(4,-3){0.08}}
\put(0.675,3.2){$\alpha$}
\put(-0.2,4.5){\vector(3,4){0.61}}
\put(1.3,6.4){\vector(-3,-4){0.55}}
\put(1.1,5.75){$t=0.06\,\mathrm{m}$}
\put(2.126,4.332){\line(0,-1){4}}
\put(0,0.75){\vector(1,0){2.126}}
\put(0.75,1){$\rho_0$}
\end{picture}
\caption{The Girkmann problem. Cross-section of the structure.} 
\label{fig:problem}
\end{center}
\end{figure}
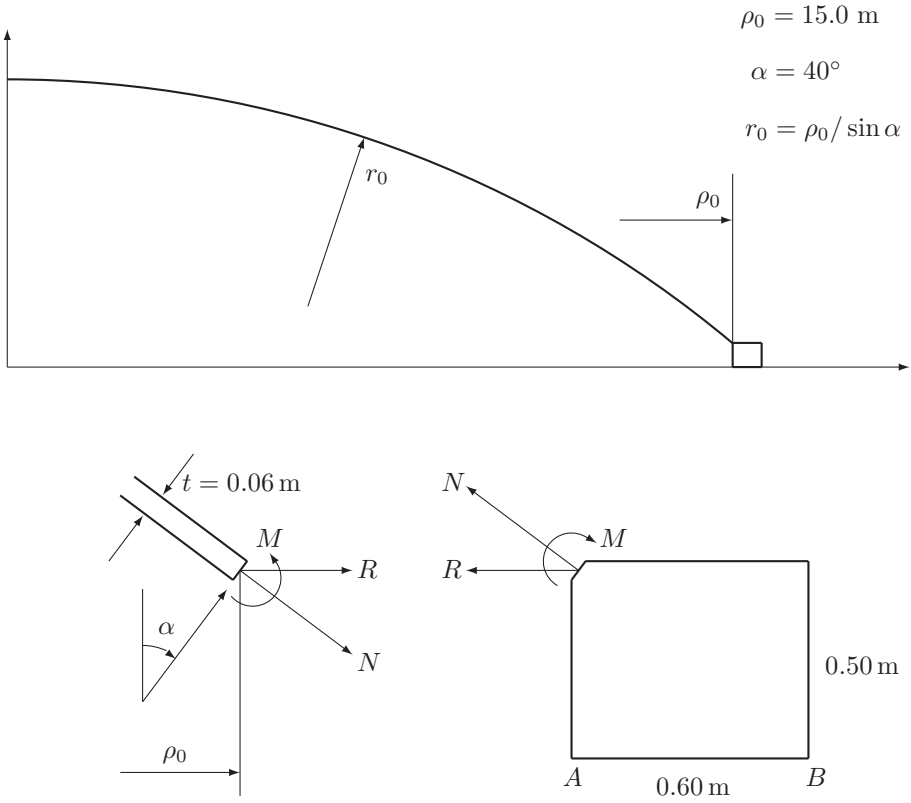

We start by recalling the statement of the Girkmann benchmark problem from \cite{Pitkaranta2008,Szabo2009}. The problem involves a concrete structure consisting of a spherical dome stiffened by a foot ring under a dead gravity load, see Fig.~\ref{fig:problem}. The task is to determine the values of the transverse shear force and the meridional bending moment at the junction between the dome and the ring as well as the maximum value of the bending moment in the dome assuming that the gravity load is equilibrated by a uniform pressure acting at the base of the ring. The material of the concrete is assumed to be linearly elastic, homogeneous and isotropic with vanishing Poisson ratio. The value of the Young modulus is specified as $E=20.59\,\mathrm{GPa}$ although it has no effect on the values of the quantities of interest.

We follow here the classical approach, where the unknown reactions are taken to be the horizontal force $R$ and the bending moment $M$ and are assumed to be positive when acting on the shell. In this splitting (shown in Fig.~\ref{fig:problem}), the normal force $N$ becomes determined from the vertical force balance as
\begin{equation} \label{eqn:normal force}
N = \frac{-gr_0}{1+\cos \alpha},
\end{equation}
where $g = F t$ is the vertical surface load density corresponding to the assumed weight density $F=32 690\,\mathrm{N}/\mathrm{m}^3$, $\alpha$ is the opening angle of the dome and $r_0$ its radius (Fig.~\ref{fig:problem}). The shear force requested in the problem statement is then defined as $Q = R/\sin\alpha$. 

The classical solution procedure can then be formulated as follows. Let us assume that $\Lambda$ and $\Psi$ denote the horizontal displacement of the midpoint of the junction and the angle of rotation of the junction line, respectively. Then, the linearity of the material model implies that the following relations must hold
\begin{equation} \label{eqn:compliance constants}
\begin{aligned}
E \Lambda &= E \Lambda_0^{\text{S}} + k_{11}^{\text{S}}R + k_{12}^{\text{S}} M = E \Lambda_0^{\text{R}} + k_{11}^{\text{R}}R + k_{12}^{\text{R}} M, \\
E \Psi &= E \Psi_0^{\text{S}} + k_{21}^{\text{S}}R + k_{22}^{\text{S}} M = E \Psi_0^{\text{R}} + k_{21}^{\text{R}}R + k_{22}^{\text{R}} M,
\end{aligned}
\end{equation}
where $\Lambda_0^{\text{S/R}}$ and $\Psi_0^{\text{S/R}}$ denote the displacement and the rotation of the Shell/Ring due to known loads only and $k_{ij}^{\text{S/R}}$, $i,j=1,2$ are the compliance constants associated to the unknown reactions $R$ and $M$. These parameters can be defined separately for the shell and the ring by analysing sequentially the following load cases
\begin{align}
\tag{Case 1} &\text{Gravity load (shell), equilibrating pressure (ring) and $N$ as in \eqref{eqn:normal force},} \\
\tag{Case 2} &R=1\,\mathrm{N}/\mathrm{m}, \\
\tag{Case 3} &M=1\,\mathrm{Nm}/\mathrm{m},
\end{align} 
and recording the values of the horizontal displacement and the rotation in each case.

\subsection{Convergence studies}
\begin{figure}
\includegraphics[width=0.333\linewidth]{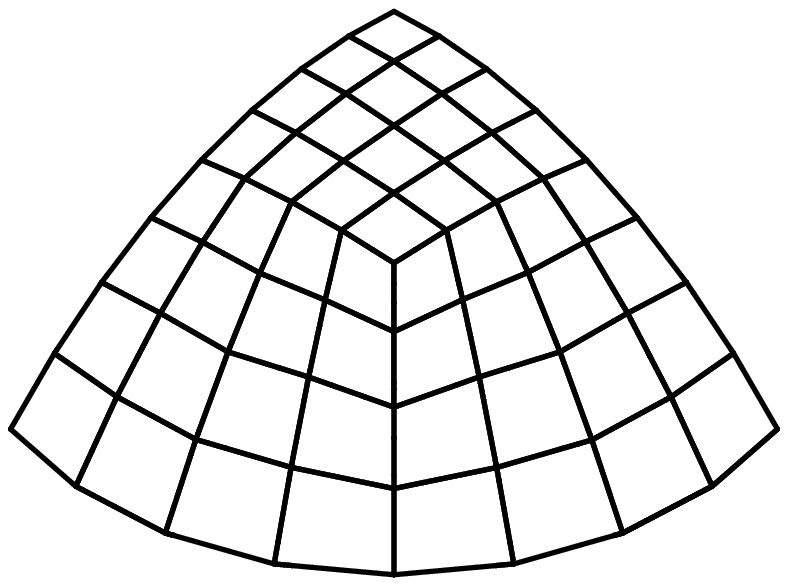}~
\includegraphics[width=0.333\linewidth]{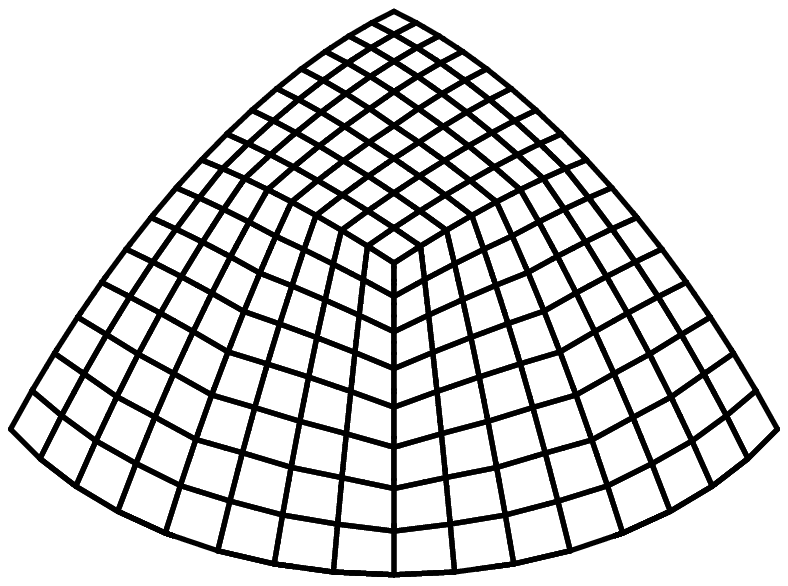}~
\includegraphics[width=0.333\linewidth]{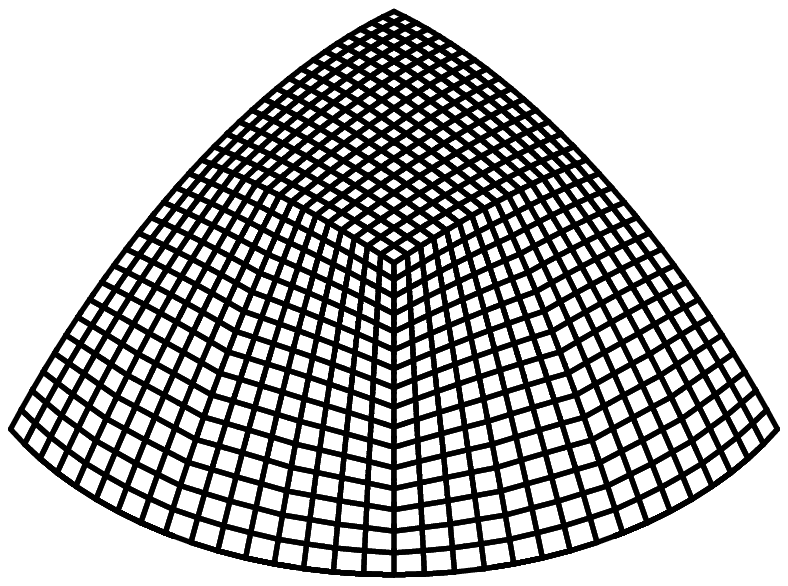}
\caption{Regularly refined mesh sequence.} \label{fig:regular_mesh}
\end{figure}

\begin{figure}
\includegraphics[width=0.333\linewidth]{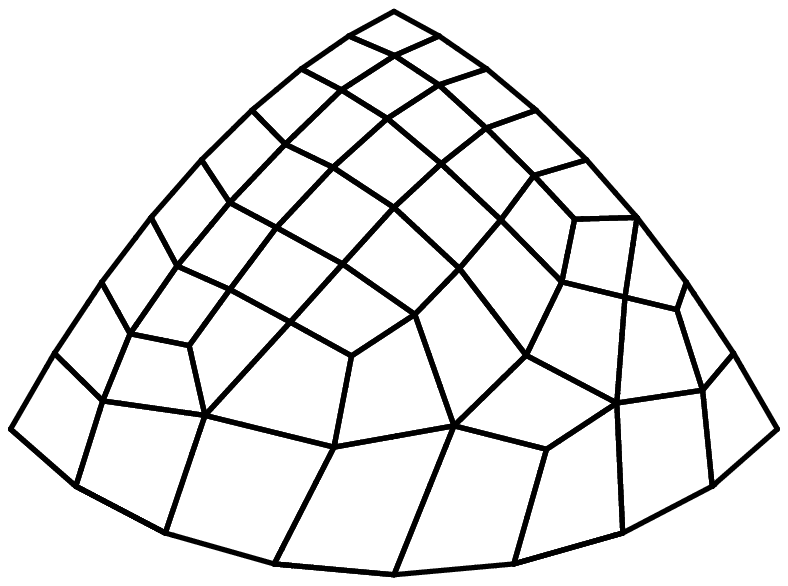}~
\includegraphics[width=0.333\linewidth]{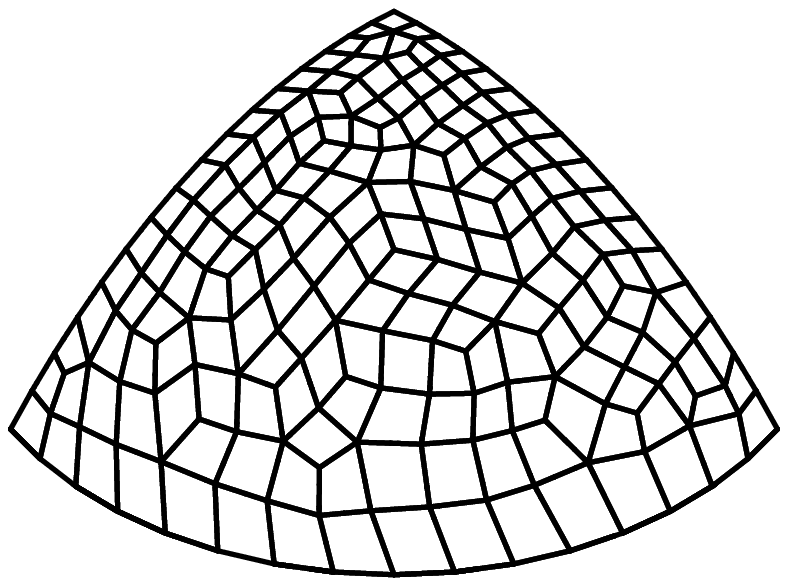}~
\includegraphics[width=0.333\linewidth]{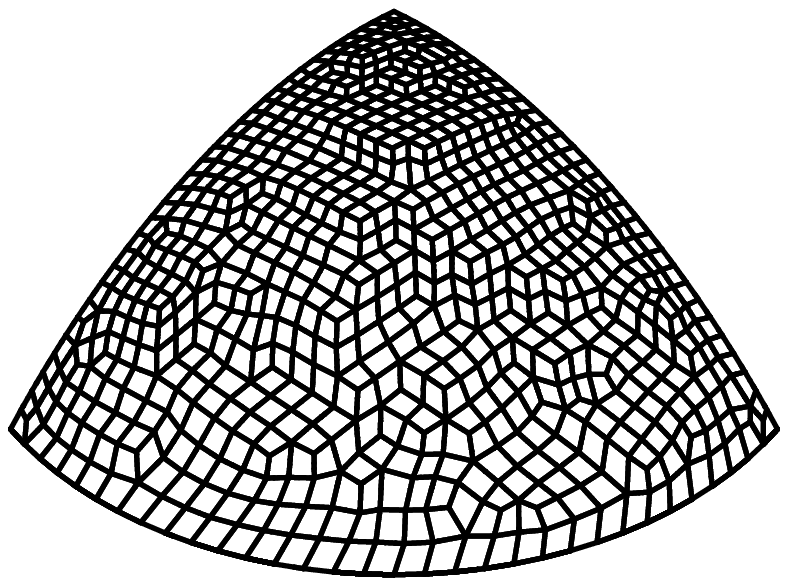}
\caption{Mesh sequence generated with the frontal method of Gmsh.} \label{fig:frontal_mesh}
\end{figure}

Our focus is on the performance of the shell elements so that we consider first the convergence of the parameters in \eqref{eqn:compliance constants} for the dome using two different kinds of mesh sequences with the maximum element size $h$ approaching zero in both cases. The first mesh sequence is based on regular refinement of the initial mesh with three elements as shown in Fig.~\ref{fig:regular_mesh}. The second mesh sequence is shown in Fig.~\ref{fig:frontal_mesh} and is generated using the frontal algorithm of Gmsh restricted so that each edge has a fixed number of elements \cite{Remacle2013}.

The results are shown in Tables \ref{tab:convergence_regular} and \ref{tab:convergence_frontal} for the different mesh sequences. The reported values are normalized against the reference values
\begin{equation} \label{eqn:compliance_coefficients_shell}
\begin{alignedat}{3}
E\Lambda_0^S &\approx \num{-2.300e6}\,\mathrm{N/m}, & \quad 
k_{11}^S &\approx \num{8.345e3} , &\quad 
k_{12}^S &\approx \num{1.477e4}\,\mathrm{1/m}, \\
E\Psi_0^S &\approx \num{-9.338e5}\,\mathrm{N/m^2}, & \quad 
k_{21}^S &\approx \num{-1.477e4}\,\mathrm{1/m} , &\quad 
k_{22}^S &\approx \num{-5.113e4}\,\mathrm{1/m^2}.
\end{alignedat}
\end{equation}
computed using the same axisymmetric shell model as in \cite{Niemi2012}. The values are nearly identical to those of \cite{Pitkaranta2012} computed using a classical shell model which neglects transverse shear deformations.\footnote{The unit of force adopted in reference \cite{Pitkaranta2012} is Girkmann's kilopond, 1 $\mathrm{G} = 9.807\,\mathrm{N}$.}

The discretization parameter $N$ stands for the number of elements per edge and is varied from 8 to 256. The required displacement and rotation values are calculated as the averages of the corresponding nodal values over the junction line and, for the stabilized variants of the methods, the value $\alpha=0.2$ is used for the stabilization parameter. 
These results do not exhibit significant differences between the formulations apart from the standard displacement method DISP4 which locks in Cases 2 and 3 as expected.
\begin{table}
\begin{center}
$E\Lambda_0^S$ \\
\begin{tabular}{lccccc}
\hline
$N$ &DISP4 		&MITC4C 	&MITC4S 	&Stab.~MITC4C 	&Stab.~MITC4C \\
\hline
8	&0.98		&0.70   	&0.75   	&0.68   	&0.75   \\
16	&1.00   	&0.89   	&0.90   	&0.89   	&0.90   \\
32	&1.00   	&0.97   	&0.97   	&0.97   	&0.97   \\
64	&1.00   	&0.99   	&0.99   	&0.99   	&0.99   \\
128	&1.00   	&1.00   	&1.00   	&1.00   	&1.00   \\
256 &1.00   	&1.00   	&1.00   	&1.00   	&1.00
\end{tabular}
\vspace{5mm}

$E\Psi_0^S$ \\
\begin{tabular}{lccccc}
\hline
$N$ &DISP4 &MITC4C &MITC4S &Stab.~MITC4C &Stab.~MITC4S \\
\hline
8	&1.00	&1.96	&2.19	&1.93	&2.12	\\
16	&0.99	&1.42	&1.47	&1.41	&1.47	\\
32	&1.00	&1.13	&1.14	&1.13	&1.14	\\
64	&1.01	&1.03	&1.04	&1.03	&1.04	\\
128	&1.01	&1.01	&1.01	&1.01	&1.01	\\
256	&1.01	&1.00	&1.00	&1.00	&1.00				
\end{tabular}
\vspace{5mm}

$k_{11}^S$ \\
\begin{tabular}{lccccc}
\hline
$N$ &DISP4 &MITC4C &MITC4S &Stab.~MITC4C &Stab.~MITC4S \\
\hline
8	&0.20	&0.51	&0.64	&0.52	&0.66	\\
16	&0.28	&0.78	&0.87	&0.78	&0.89	\\
32	&0.40	&0.93	&0.96	&0.93	&0.97	\\
64	&0.56	&0.98	&0.99	&0.98	&0.99	\\
128	&0.74	&1.00	&1.00	&1.00	&1.00	\\
256	&0.89	&1.00	&1.00	&1.00	&1.00				
\end{tabular}
\vspace{5mm}

$k_{12}^S$ \\
\begin{tabular}{lccccc}
\hline
$N$ &DISP4 &MITC4C &MITC4S &Stab.~MITC4C &Stab.~MITC4S \\
\hline
8	&0.04	&0.39	&0.56	&0.36	&0.53	\\
16	&0.09	&0.71	&0.83	&0.70	&0.83	\\
32	&0.17	&0.91	&0.95	&0.90	&0.95	\\
64	&0.32	&0.97	&0.99	&0.97	&0.99	\\
128	&0.56	&0.99	&1.00	&0.99	&1.00	\\
256	&0.80	&1.00	&1.00	&1.00	&1.00				
\end{tabular}
\vspace{5mm}

$k_{22}^S$ \\
\begin{tabular}{lccccc}
\hline
$N$ &DISP4 &MITC4C &MITC4S &Stab.~MITC4C &Stab.~MITC4S \\
\hline
8	&0.01	&0.52	&0.66	&0.67	&0.79	\\
16	&0.03	&0.78	&0.88	&0.81	&0.90	\\
32	&0.07	&0.93	&0.97	&0.94	&0.97	\\
64	&0.18	&0.98	&0.99	&0.98	&0.99	\\
128	&0.42	&1.00	&1.00	&1.00	&1.00	\\
256	&0.72	&1.00	&1.00	&1.00	&1.00				
\end{tabular}

\caption{Convergence of the compliance coefficients for the dome with respect to the mesh parameter $N$ on the \emph{regular mesh sequence (Fig.~\ref{fig:regular_mesh})}.}
\label{tab:convergence_regular}
\end{center}
\end{table}
\begin{table}
\begin{center}
$E\Lambda_0^S$ \\
\begin{tabular}{lccccc}
\hline
$N$ &DISP4 		&MITC4C 	&MITC4S 	&Stab.~MITC4C 	&Stab.~MITC4C \\
\hline
8	&1.00   	&0.74   	&0.92   	&0.72   	&0.98	\\
16	&1.02   	&0.91   	&0.94   	&0.90   	&0.94	\\
32	&1.01   	&0.97   	&0.98   	&0.97   	&0.98  	\\
64	&1.01   	&0.99   	&1.00   	&0.99   	&1.00  	\\
128	&1.00   	&1.00   	&1.00   	&1.00   	&1.00  	\\
256	&1.00   	&1.00   	&1.00   	&1.00   	&1.00   
\end{tabular}
\vspace{5mm}

$E\Psi_0^S$ \\
\begin{tabular}{lccccc}
\hline
$N$ &DISP4 &MITC4C &MITC4S &Stab.~MITC4C &Stab.~MITC4S \\
\hline
8	&0.96	&1.78	&1.83	&1.83	&1.81	\\
16	&0.96	&1.37	&1.45	&1.37	&1.45	\\
32	&0.98	&1.12	&1.05	&1.12	&1.05	\\
64	&0.99	&1.03	&1.01	&1.03	&1.01	\\
128	&1.00	&1.01	&1.00	&1.01	&1.00	\\
256	&1.00	&1.00	&1.00	&1.00	&1.00		
\end{tabular}
\vspace{5mm}

$k_{11}^S$ \\
\begin{tabular}{lccccc}
\hline
$N$ &DISP4 &MITC4C &MITC4S &Stab.~MITC4C &Stab.~MITC4S \\
\hline
8	&0.20	&0.48	&0.57	&0.50	&0.68	\\
16	&0.28	&0.71	&0.78	&0.75	&0.82	\\
32	&0.40	&0.89	&0.91	&0.91	&0.94	\\
64	&0.54	&0.96	&0.98	&0.97	&0.98	\\
128	&0.73	&0.99	&0.99	&0.99	&1.00	\\	
256	&0.86	&1.00	&1.00	&1.00	&1.00	
\end{tabular}
\vspace{5mm}

$k_{12}^S$ \\
\begin{tabular}{lccccc}
\hline
$N$ &DISP4 &MITC4C &MITC4S &Stab.~MITC4C &Stab.~MITC4S \\
\hline
8	&0.04	&0.32	&0.42	&0.33	&0.46	\\
16	&0.08	&0.62	&0.70	&0.63	&0.72	\\
32	&0.17	&0.85	&0.88	&0.87	&0.90	\\
64	&0.30	&0.95	&0.97	&0.96	&0.97	\\
128	&0.53	&0.99	&0.99	&0.99	&0.99	\\
256	&0.75	&1.00	&1.00	&1.00	&1.00		
\end{tabular}
\vspace{5mm}

$k_{22}^S$ \\
\begin{tabular}{lccccc}
\hline
$N$ &DISP4 &MITC4C &MITC4S &Stab.~MITC4C &Stab.~MITC4S \\
\hline
8	&0.01	&0.34	&0.41	&0.69	&0.78	\\
16	&0.02	&0.70	&0.76	&0.77	&0.84	\\
32	&0.06	&0.88	&0.90	&0.91	&0.94	\\
64	&0.16	&0.96	&0.97	&0.97	&0.98	\\
128	&0.37	&0.99	&0.99	&0.99	&0.99	\\
256	&0.64	&1.00	&1.00	&1.00	&1.00		
\end{tabular}

\caption{Convergence of the compliance coefficients for the dome with respect to the mesh parameter $N$ on the \emph{frontal mesh sequence (Fig.~\ref{fig:frontal_mesh})}.}
\label{tab:convergence_frontal}
\end{center}
\end{table}

\begin{figure}
\includegraphics[width=0.5\linewidth]{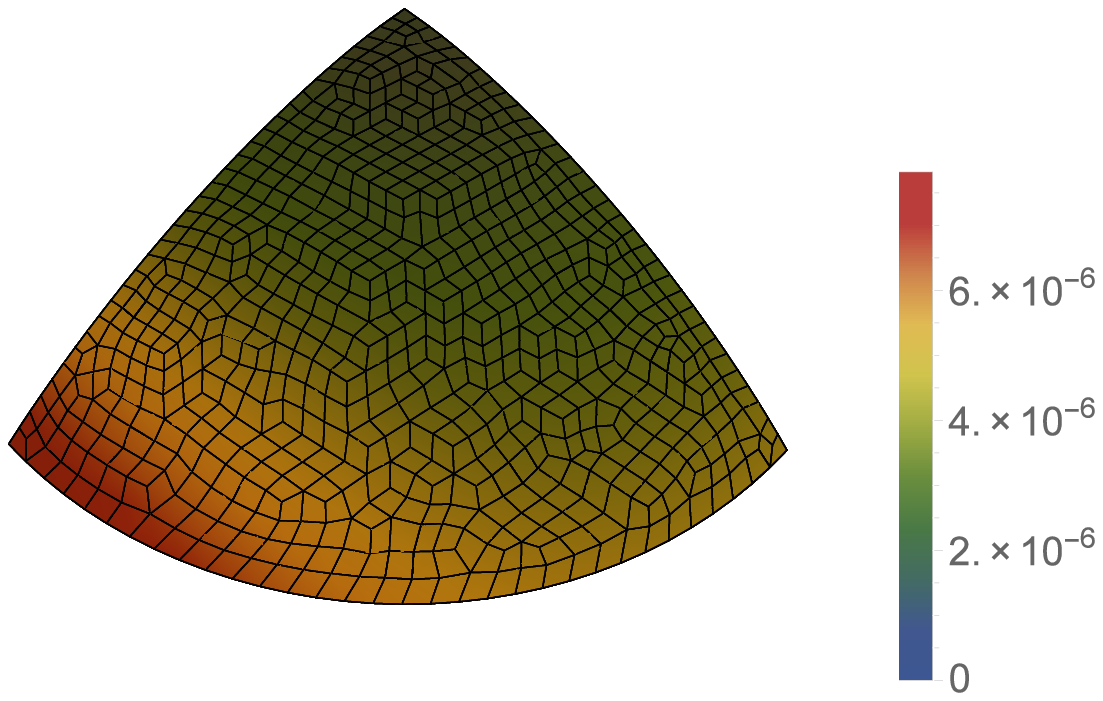}~
\includegraphics[width=0.5\linewidth]{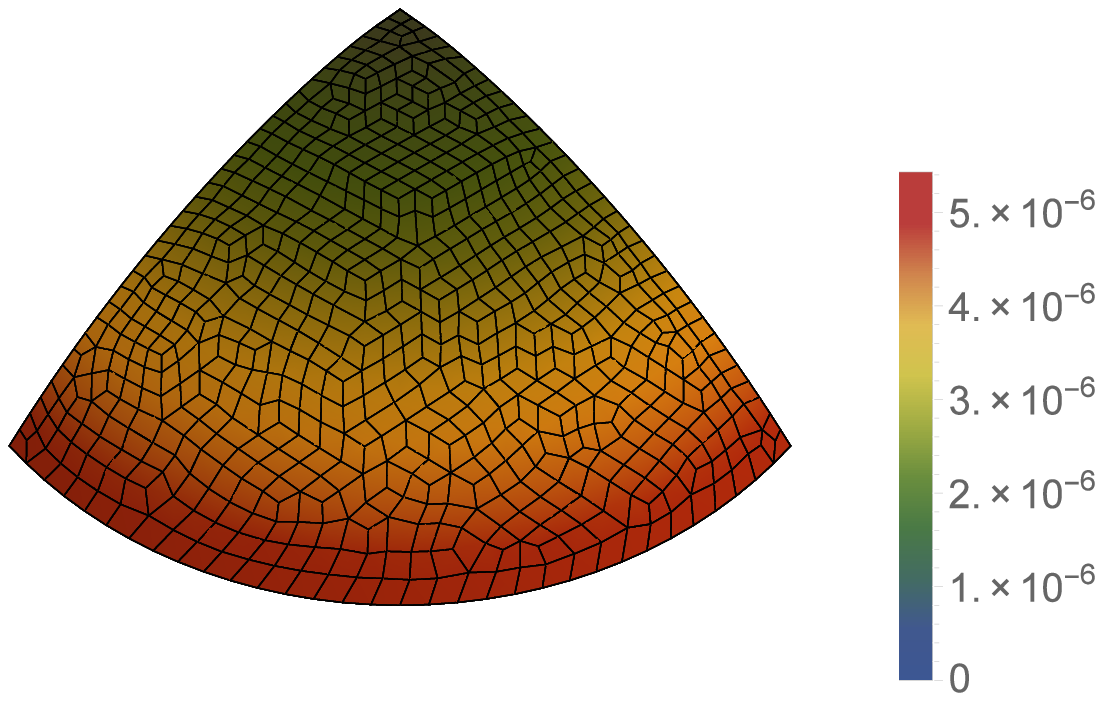}
\caption{Total displacement in meters in the membrane-dominated load Case 1. MITC4S (left) versus MITC4C (right) on the frontal mesh with $N=32$.}
\label{fig:Case_1_total_diplacement_frontal}
\end{figure}

\begin{figure}
\includegraphics[width=0.5\linewidth]{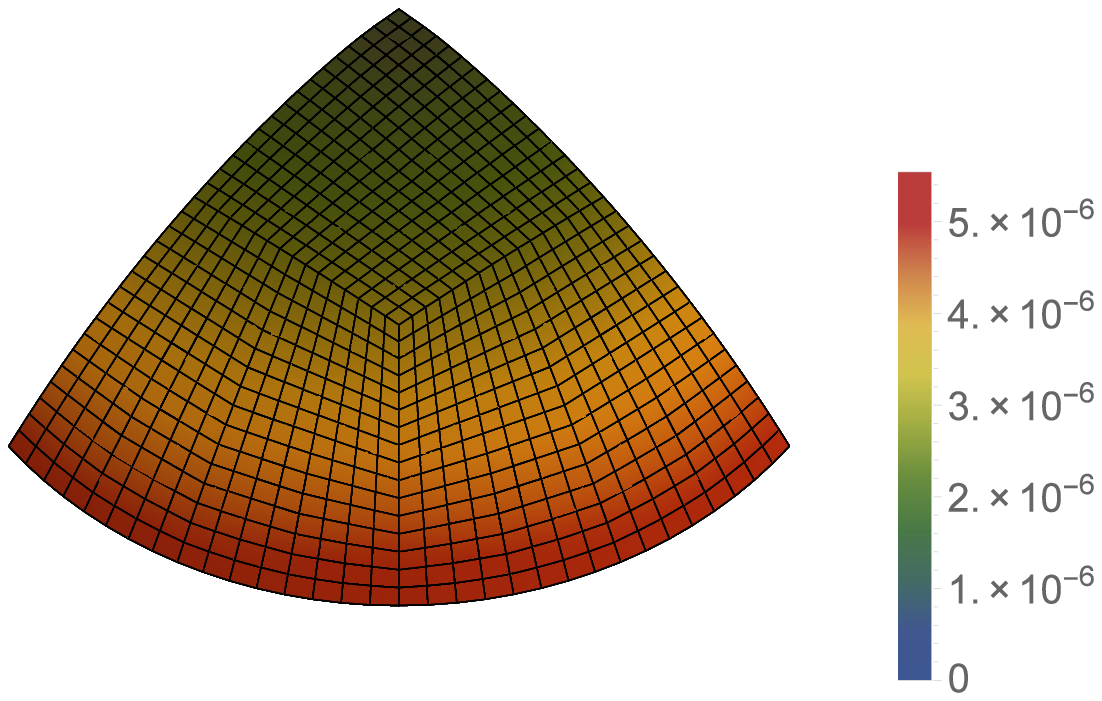}~
\includegraphics[width=0.5\linewidth]{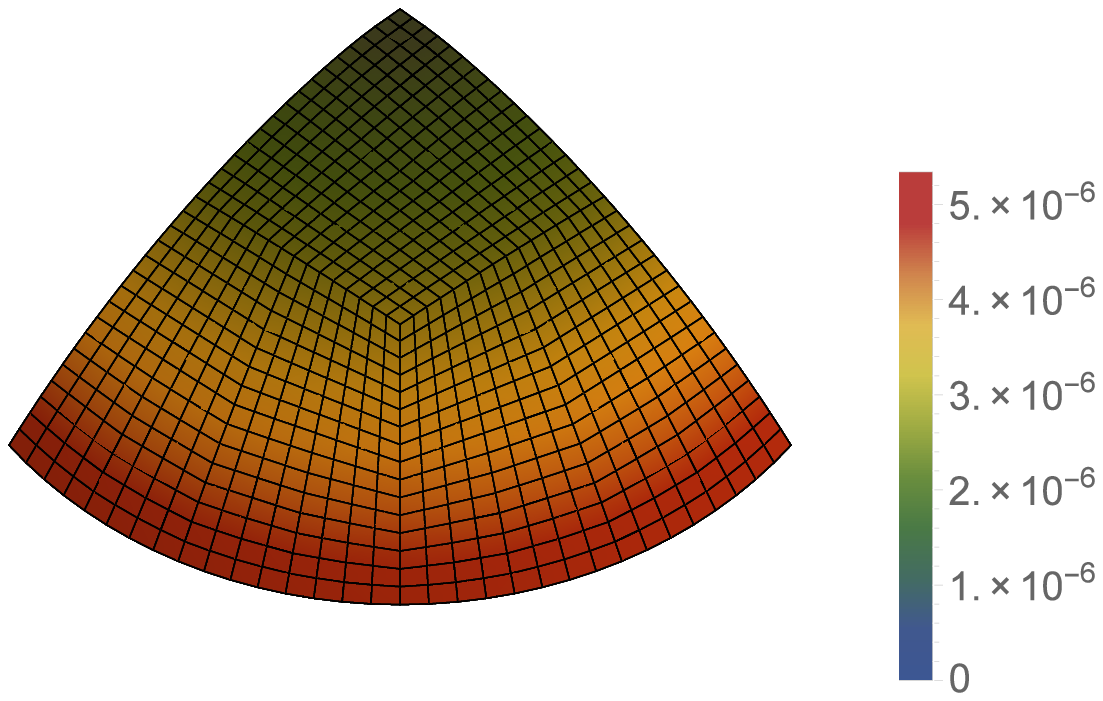}
\caption{Total displacement in meters in the membrane-dominated load Case 1. MITC4S (left) versus MITC4C (right) on the regular mesh with $N=32$.}
\label{fig:Case_1_total_diplacement_regular}
\end{figure}

However, a more detailed investigation reveals a rather drastic difference between the different formulations in the solutions of the membrane-dominated Case 1. Fig.~\ref{fig:Case_1_total_diplacement_frontal} shows the total displacement of the shell mid-surface calculated on the frontal mesh with 16 elements per edge using the MITC4S and MITC4C formulations. It is evident that the MITC4S formulation suffers from numerical instabilities associated to the consistency error arising from the reduction of the membrane strains, cf.~\cite{Niemi2008a}. The instability is manifested here by the loss of symmetry of the numerical solution. It is intriguing that the circumferentially averaged displacements still converge and that the phenomenon disappears when a regular mesh is used as shown in Fig.~\ref{fig:Case_1_total_diplacement_regular}.

Finally, it should be pointed out that similar instabilities as shown in Fig.~4 (left) are featured also by the lowest-order linear or bilinear shell elements employed currently in many industrial FEA programs.

\subsection{Section force and moment computations}
In order to determine the unknown reactions $R$ and $M$ and the deformation of the stiffened shell, we need the compliance coefficients of the ring. These can be approximated by using the principle of virtual work as
\begin{equation} \label{eqn:compliance_coefficients_ring}
\begin{alignedat}{3}
E\Lambda_0^R &\approx \num{1.363e7}\,\mathrm{N/m}, & \quad 
k_{11}^R &\approx -2683, &\quad 
k_{12}^R &\approx 8418\,\mathrm{1/m},\\
E\Psi_0^R &\approx \num{-6.949e6}\,\mathrm{N/m^2}, & \quad 
k_{21}^R &\approx \num{-8418}\,\mathrm{1/m}, &\quad 
k_{22}^R &\approx \num{3.696e4}\,\mathrm{1/m^2}.
\end{alignedat}
\end{equation}
These values are obtained by taking the aforementioned displacements $\Lambda$ and $\Psi$ as the only degrees of freedom as in \cite{Pitkaranta2012}. This corresponds to the kinematic assumption that the ring cross section deforms as a rigid body. The associated $2 \times 2$ stiffness matrix is evaluated by using numerical integration up to the machine precision in cylindrical coordinates over the exact pentagonal shape of the ring. The ring is assumed weightless here as in the original treatment by Girkmann and in the contemporary verification challenge.

Substitution of \eqref{eqn:compliance_coefficients_ring} and \eqref{eqn:compliance_coefficients_shell} into \eqref{eqn:compliance constants} yields the values of the unknown horizontal section force and bending moment:
\begin{equation} \label{eqn:unknown reactions}
R \approx 1467\,\mathrm{N/m} \quad \& \quad M \approx -37.36\,\mathrm{Nm/m}.
\end{equation}
Finally, the deformation of the stiffened dome can be computed by solving the shell problem with the active loads as in
\begin{equation} 
\tag{Case 4} \text{Gravity load, $N$ according to \eqref{eqn:normal force} and $R,M$ according to \eqref{eqn:unknown reactions}}.
\end{equation}

\begin{figure}
\includegraphics[width=0.5\linewidth]{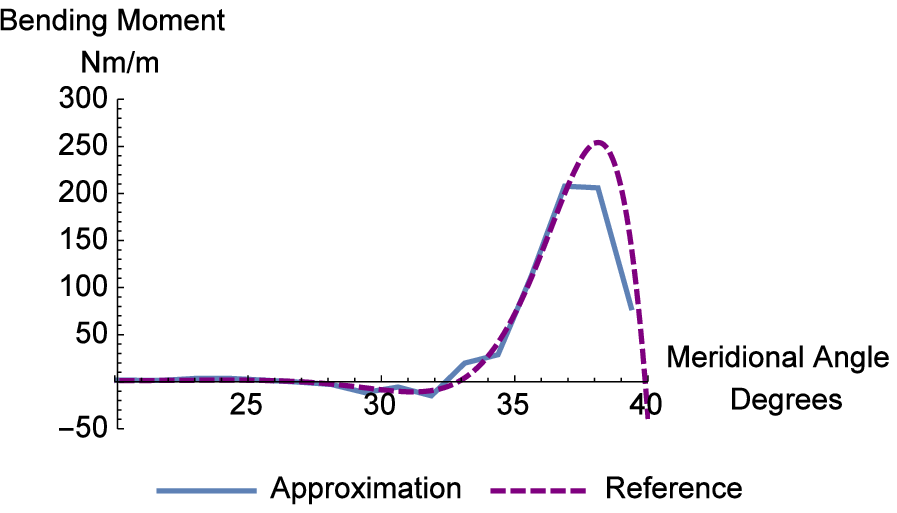}~
\includegraphics[width=0.5\linewidth]{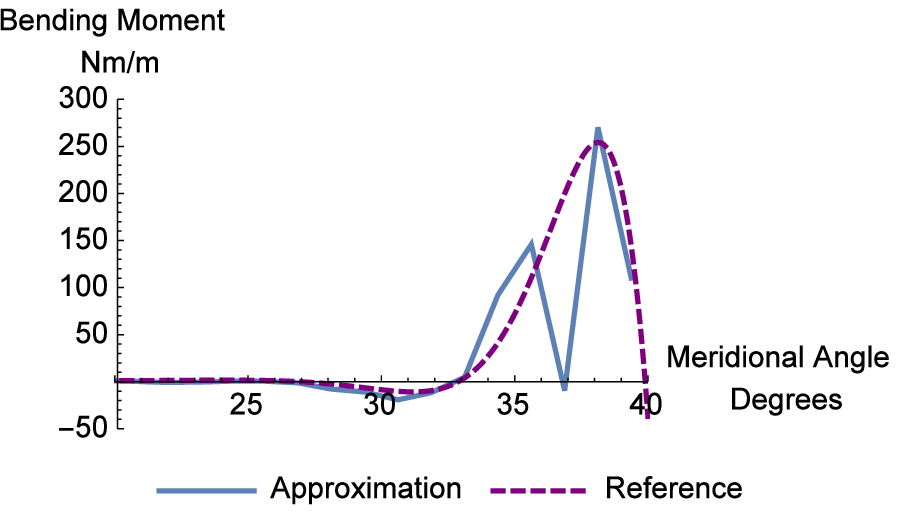}
\caption{Distribution of the meridional bending moment in the stiffened dome calculated along the left and right edges of the computational domain, respectively. MITC4S formulation on the \emph{frontal mesh} with $N=32$.}
\label{fig:Moment_MITC4S_32GQ}
\end{figure}

\begin{figure}
\includegraphics[width=0.5\linewidth]{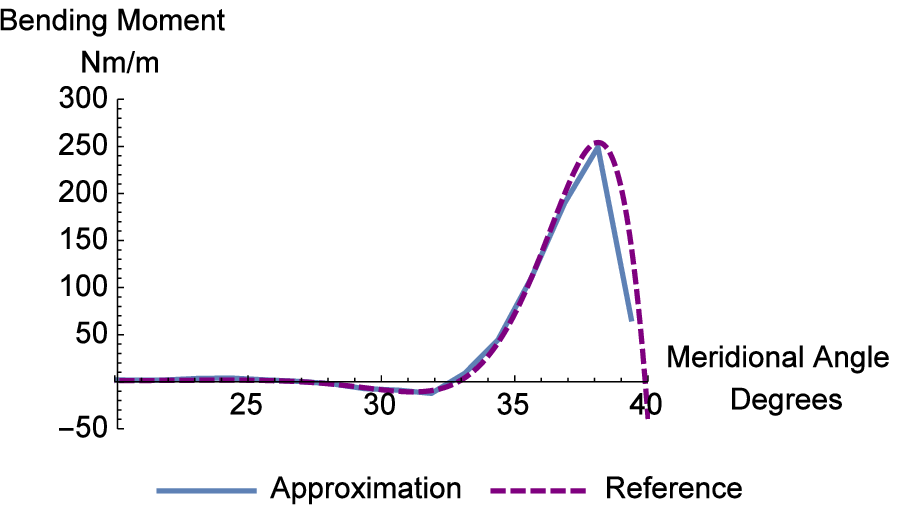}~
\includegraphics[width=0.5\linewidth]{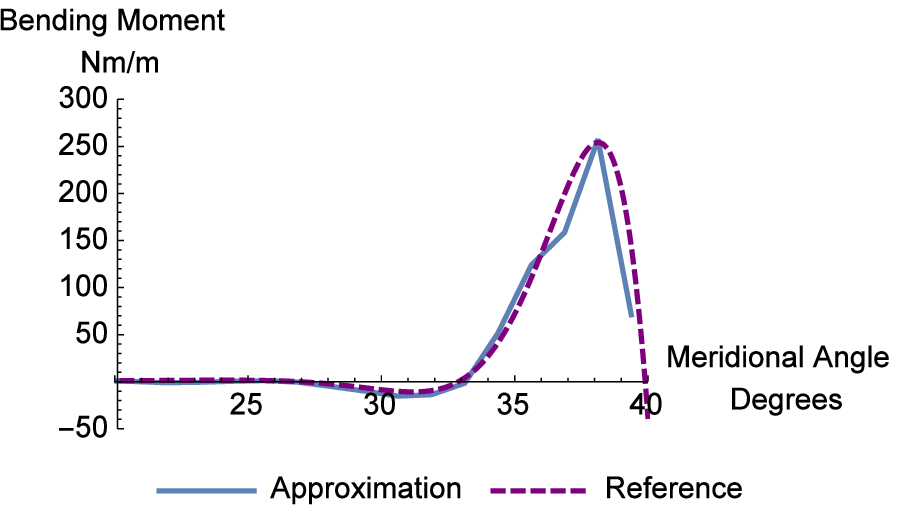}
\caption{Distribution of the meridional bending moment in the stiffened dome calculated along the left and right edges of the computational domain, respectively. \emph{Stabilized} MITC4S formulation on the \emph{frontal mesh} with $N=32$.}
\label{fig:Moment_Stab_MITC4S_32GQ}
\end{figure}

In order to demonstrate the influence of the stabilization technique \eqref{eqn:shear balancing}, we show in Figs.~\ref{fig:Moment_MITC4S_32GQ}--\ref{fig:Moment_Stab_MITC4S_32GQ} the distribution of the meridional bending moment as computed with the different formulations along the left and right edges of the computational domain. The post-processing is carried out directly from the nodal rotations and only the range $(20^\circ,40^\circ)$ is shown here since the bending effects are confined to a narrow region near the edge.

\begin{figure}
\includegraphics[width=0.5\linewidth]{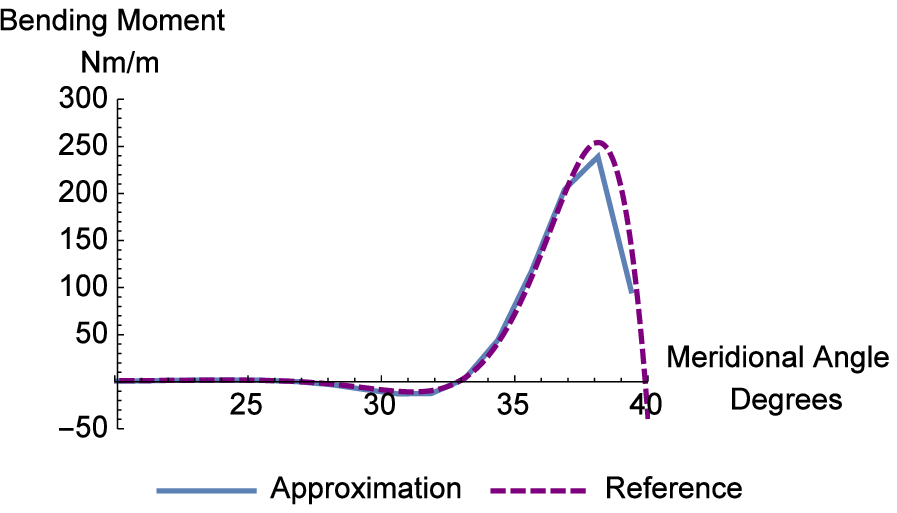}~
\includegraphics[width=0.5\linewidth]{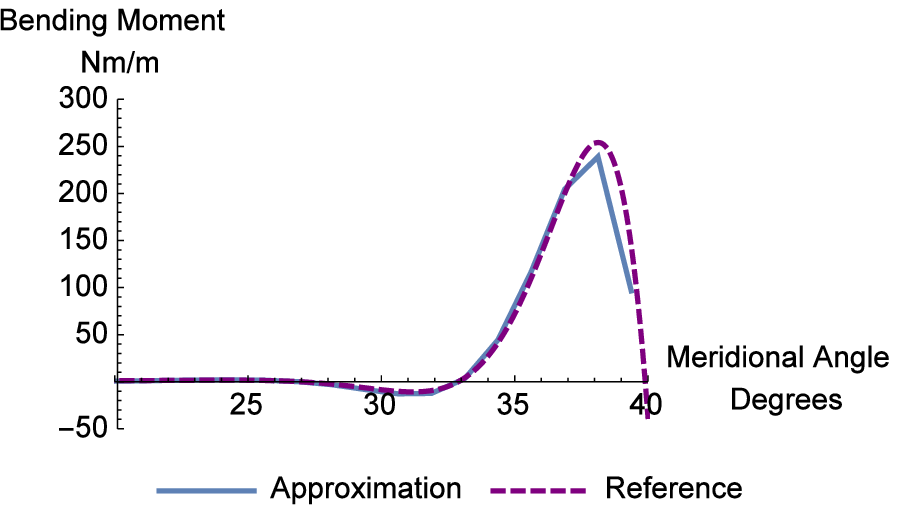}
\caption{Distribution of the meridional bending moment in the stiffened dome calculated along the left and right edges of the computational domain, respectively. MITC4S formulation on the \emph{regular mesh} with $N=32$.}
\label{fig:Moment_MITC4S_32}
\end{figure}

In this case, there is no visible difference between the MITC4S and MITC4C formulations. For instance, on the frontal mesh with $N=32$, both formulations feature non-physical oscillations but the stabilization technique \eqref{eqn:shear balancing} improves the results as shown in Figs.~\ref{fig:Moment_MITC4S_32GQ} and \ref{fig:Moment_Stab_MITC4S_32GQ}. On the other hand, a feasible solution is obtained even without the stabilization when the mesh is regular as shown in Fig.~\ref{fig:Moment_MITC4S_32}.

\section{Conclusions} \label{sec:conclusions}
We have presented a finite element framework for analysing the structural response of thin elastic shells. The element stiffness matrices are computed according to shell theory taking into account transverse shear deformations. The framework enables enables explicit reduction of the membrane strains and the transverse shear strains in order to resolve numerical locking problems.

Different variants of quadrilateral elements have been examined in the Girkmann problem and a detailed verification study has been performed. We have demonstrated that the MITC4S element with reduced membrane strains may feature numerical instabilities in membrane-dominated situations if the mesh is irregular. The MITC4C where only the transverse shear strains are reduced, does not have this drawback. On the other hand, the bending moment in the Girkmann problem can be calculated accurately with both methods on regular meshes. On irregular meshes some non-physical oscillations occur but these can be avoided by employing  the transverse shear balancing as a stabilization technique.

\bibliographystyle{siam}
\bibliography{library}
\end{document}